\documentclass[12pt,a4paper]{amsart}
\usepackage{amsfonts}
\usepackage{amscd}
\usepackage{amssymb}
\usepackage{mathrsfs}
\usepackage[scale=0.75]{geometry}

\newtheorem{thm}{Theorem}

\theoremstyle{definition}
\newtheorem{defn}[thm]{Definition}

\theoremstyle{remark}
\newtheorem{rem}[thm]{Remark}

\newcommand{\C}{\mathbb C}

\newcommand{\HC}{\mathbb H}

\author{Anton Mellit}
\title{Angle structures on $3$-manifolds}
\email{anton.mellit@univie.ac.at}
\address{Faculty of Mathematics, University of Vienna, \\
	Oskar-Morgenstern-Platz 1, 1090 Vienna, Austria}
\date{\today}
\thanks{Partially supported by the projects Y963-N35 and
	P-31705 of the Austrian Science Fund. This note was written in 2007 during the author's stay at MPIfM. It is uploaded to arXiv at the request of Herbert Gangl.}

\begin{document}
\begin{abstract}
	Given a compact oriented triangulated $3$-manifold we find a non-trivial condition satisfied by certain labelings of the tetrahedra by elements of an arbitrary abelian group which we call angle structures. Smoothness of the manifold is used in an essential way. This is inspired by the notion of the volume of hyperbolic manifolds, which would correspond to the case when the abelian group is the multiplicative group of $\mathbb{C}$, but the construction here seems to be more general, in particular it only uses the abelian group structure.
\end{abstract}
\bibliographystyle{alpha}
\maketitle

Let $A$ be an abelian group. Let $X$ be a trangulated oriented $3$-dimensional manifold. Let $\Delta=(a,b,c,d)$ be an oriented $3$-dimensional simplex of the triangulation.
\begin{defn}
An angle structure on $\Delta$ a collection of $12$ elements of $A$, denoted $k_{xy}$, one for each oriented edge $xy$ satisfying relations:
\[
k_{xy}=k_{yx},\; \sum_{y} k_{xy}=0 \;\text{for any vertex $x$.}
\]
\end{defn}

We have
\[
k_{ab}+k_{ac}+k_{ad}=0,\; k_{ab}+k_{bc}+k_{bd}=0,\; k_{ac}+k_{bc}+k_{cd}=0,\; k_{ad}+k_{bd}+k_{cd}=0.
\]
This implies
\[
2(k_{ab}-k_{cd})=2(k_{ac}-k_{bd})=2(k_{ad}-k_{bc})=0.
\]

Suppose each tetrahedron of the triangulation has an angle structure. The corresponding elements will be denoted $k_{ab}^{cd}\in A$, i.e. for a simplex $\Delta=(a,b,c,d)$ we have $12$ elements listed below:
\[
k_{ab}^{cd}=k_{ba}^{dc}, k_{ac}^{db}=k_{ca}^{bd}, k_{ad}^{bc}=k_{da}^{cb}, k_{bc}^{ad}=k_{cb}^{da}, k_{bd}^{ca}=k_{db}^{ac}, k_{cd}^{ab}=k_{dc}^{ba}.
\]
Therefore the notation is invariant with respect to even permutations of vertices of $\Delta$. 

\begin{defn}
For each edge $ab$ of $X$ let $(ab x_i x_{i+1})$ be the $3$-simplices adjacent to $ab$ with $i=0,1,\ldots,n-1$, and $x_n=x_0$. Suppose the following condition is satisfied:
\[
\sum_i k_{ab}^{x_i x_{i+1}}=0.
\]
Then we call our manifold \emph{angled}.
\end{defn}

Suppose we have a finite set of tetrahedra with angle structure. Is it possible to glue them together and obtain an angled oriented manifold? We provide a necessary condition.
\begin{thm}
In an angled oriented manifold the following condition is satisfied in the group $\Lambda^2 A$:
\[
\sum_{abcd} 4 k_{ab}^{cd}\wedge k_{ac}^{db}+k_{ad}^{bc}\wedge k_{ad}^{bc}+k_{bc}^{ad}\wedge k_{bc}^{ad}=0,
\]
the sum is over $3$-simplices which compose the fundamental class of the manifold.
\end{thm}
\begin{rem}\label{rem:bloch}
Let $A=\C^\times$ and each tetrahedron is realized as an ideal tetrahedron in $\HC^3$ with cross-ratio $z$. Then its angles are $z$, $\frac{1}{1-z}$, $1-\frac{1}{z}$. To make the product $1$ one can change angles to $-z$, $\frac{1}{z-1}$, $\frac{1}{z}-1$. Then we see that up to $2$-torsion the sum of tensors $z\wedge (1-z)$ is $0$, which is well known. In other words, hyperbolic $3$-manifolds provide elements in the Bloch group, see \cite{NZ}, \cite{DS}. However our approach seems to be more general.
\end{rem}
The rest of this paper provides a proof of the theorem.

Suppose $X$ is angled. Then we can construct elements $h_{ab}^c\in A$ with the property
\[
k_{ab}^{cd} = h_{ab}^d-h_{ab}^c
\]
for each oriented simplex $(a,b,c,d)$. If $h'^c_{ab}$ is another such family then there is a family $q_{ab}$ with 
\[
h'^c_{ab}-h_{ab}^c = q_{ab}.
\]

Fix a vertex $a$. Let $ab$ and $ab'$ be edges. Join $ab$ and $ab'$ by a sequence of triangles $(a b_i b_{i+1})$ for $i=0,1,\ldots,n-1$, $b_0=b$, $b_n=b'$. Put 
\[
m_a(b,b') = \sum_{i} (h_{a b_i}^{b_{i+1}} - h_{a b_{i+1}}^{b_i}).
\]
This does not depend on the choice of the sequence since for any oriented $3$-simplex $(a b c d)$ we have
\[
(h_{a b}^{c} - h_{a c}^b) + (h_{ac}^d-h_{ad}^c) + (h_{ad}^b-h_{ab}^d)=
-k_{ab}^{cd}-k_{ac}^{db}-k_{ad}^{bc}=0.
\]
Therefore there is a family $q_{ab}\in A$ with the property
\[
m_a(b,b')=q_{ab'}-q_{ab}.
\]
If $q'_{ab}$ is another such family there exists a family $p_a\in A$ with
\[
q'_{ab}-q_{ab}=p_a.
\]
In particular for any $2$-simplex $(abc)$ we have
\[
h_{a b}^c-h_{ac}^b = q_{ac}-q_{ab}.
\]
We see that we can replace $h_{a b}$ with $h_{a b}^c+q_{ab}$ to make $h$ satisfying
\[
h_{ab}^c=h_{ac}^b\;\text{for any $2$-simplex $(abc)$.}
\]
If $h'^c_{ab}$ is another family with this condition then there is a family $p_a\in A$ with the property
\[
h'^c_{ab}-h^c_{ab}=p_a.
\]
Put $h'^c_{ab}=-h^c_{ba}$. Then for any oriented $3$-simplex $(abcd)$
\[
k_{ab}^{cd} = k_{ba}^{dc}=h_{ba}^c-h_{ba}^d=h'^d_{ab}-h'^c_{ab}.
\]
Therefore there exists a family $q_{ab}$ (this $q$ is different from the one used before) with property
\[
h^c_{ab}-h'^c_{ab} = q_{ab}.
\]
This means that for any $2$-simplex $(abc)$ we have
\[
h^c_{ab}+h^c_{ba} = q_{ab}.
\]
Let us summarize the properties of $h^c_{ab}$:
\[
h_{ab}^c=h_{ac}^b,\;h^c_{ab}+h^c_{ba} = q_{ab},\;k_{ab}^{cd} = h_{ab}^d-h_{ab}^c.
\]
Let $(abc)$ be a $2$-simplex. Then
\[
q_{ab}+q_{ac}-q_{bc}=h_{ab}^c+h_{ba}^c+h_{ac}^b+h_{ca}^b-h_{bc}^a-h_{cb}^a = 2 h_{ab}^c.
\]
Consider the following element in $\Lambda^2 A$:
\[
\phi_{abc}=h_{ab}^c\wedge h_{bc}^a + h_{bc}^a\wedge h_{ca}^b + h_{ca}^b\wedge h_{ab}^c.
\]
This element is invariant under cyclic permutations:
\[
\phi_{abc}=\phi_{bca}=\phi_{cab}.
\]
Moreover,
\[
\phi_{acb}=h_{ac}^b\wedge h_{cb}^a + h_{cb}^a\wedge h_{ba}^c + h_{ba}^c\wedge h_{ac}^b=h_{ab}^c\wedge h_{ca}^b + h_{ca}^b\wedge h_{bc}^a + h_{bc}^a\wedge h_{ab}^c=-\phi_{abc}.
\]
Since $h_{bc}^a=q_{ab}-h_{ab}^c$ and $h_{ca}^b=q_{ac}-h_{ab}^c$,
we can also write $\phi_{abc}$ as
\begin{multline*}
\phi_{abc}=2 h_{ab}^c\wedge(q_{ab}-q_{ac}) +q_{ab}\wedge q_{ac} - h_{ab}^c\wedge h_{ab}^c \\
=q_{ac}\wedge q_{ab} +q_{ab}\wedge q_{bc} +q_{bc}\wedge q_{ac} + q_{ab}\wedge q_{ab} + q_{ac}\wedge q_{ac} + h_{ab}^c\wedge h_{ab}^c\\
=q_{ac}\wedge q_{ab} +q_{ab}\wedge q_{bc} +q_{bc}\wedge q_{ac} + q_{bc}\wedge q_{bc} + h_{ab}^c\wedge h_{ab}^c.
\end{multline*}

Let 
\[
\phi_{abc}^0=q_{ac}\wedge q_{ab} +q_{ab}\wedge q_{bc} +q_{bc}\wedge q_{ac},\; \phi_{abc}^1=q_{bc}\wedge q_{bc} + h_{ab}^c\wedge h_{ab}^c.
\]

For any oriented $3$-simplex $(abcd)$ put
\[
\psi^*_{abcd}=\phi^*_{bcd}-\phi^*_{acd}+\phi^*_{abd}-\phi^*_{abc}.
\]
Then
\begin{multline*}
\psi^0_{abcd}=q_{bc}\wedge q_{cd}+q_{cd}\wedge q_{bd}+q_{bd}\wedge q_{bc}-
q_{ac}\wedge q_{cd}-q_{cd}\wedge q_{ad}-q_{ad}\wedge q_{ac}\\
+q_{ab}\wedge q_{bd}+q_{bd}\wedge q_{ad}+q_{ad}\wedge q_{ab}-
q_{ab}\wedge q_{bc}-q_{bc}\wedge q_{ac}-q_{ac}\wedge q_{ab}\\
=q_{ab}\wedge(q_{bd}-q_{ad}-q_{bc}+q_{ac})
+q_{cd}\wedge(-q_{bc}+q_{bd}+q_{ac}-q_{ad})
+(q_{bd}+q_{ac})\wedge(q_{bc}+q_{ad})\\
=(q_{bc}+q_{ad})\wedge(q_{ab}+q_{cd}) + (q_{ab}+q_{cd})\wedge(q_{bd}+q_{ac})+(q_{bd}+q_{ac})\wedge(q_{bc}+q_{ad}).
\end{multline*}

We may rewrite
\begin{align*}
q_{bc}+q_{ad} &=h_{bc}^a+h_{cb}^a+h_{ad}^c+h_{da}^c = h_{bc}^a+h_{dc}^a + q_{ac}-k_{ca}^{bd},\\
q_{ab}+q_{cd} &=h_{ab}^c+h_{ba}^c+h_{cd}^a+h_{dc}^a = h_{bc}^a+h_{dc}^a + q_{ac}+k_{ca}^{bd}.
\end{align*}
Therefore
\begin{multline*}
\psi^0_{abcd}=2 k_{ca}^{bd}\wedge(q_{bd}+q_{ac})-2 k_{ca}^{bd}\wedge(h_{bc}^a+h_{dc}^a + q_{ac})+k_{ca}^{bd}\wedge k_{ca}^{bd}+h_{bc}^a\wedge h_{bc}^a+h_{dc}^a\wedge h_{dc}^a+q_{ac}\wedge q_{ac}\\
=2k_{ca}^{bd}\wedge(q_{bd}-h_{bc}^a-h_{dc}^a) + q_{bd}\wedge q_{bd} + q_{ac}\wedge q_{ac}.
\end{multline*}
Take into account that
\[
q_{bd}-h_{bc}^a-h_{dc}^a = h_{bd}^a+h_{db}^a-h_{bc}^a-h_{dc}^a=h_{ba}^d+h_{da}^b-h_{ba}^c-h_{da}^c=-k_{ba}^{dc}+k_{da}^{cb},
\]
we obtain
\[
\psi^0_{abcd}=2 k_{ca}^{bd}\wedge(k_{da}^{cb}-k_{ba}^{dc})+q_{bd}\wedge q_{bd} + q_{ac}\wedge q_{ac}=4 k_{ab}^{cd}\wedge k_{ac}^{db}+q_{bd}\wedge q_{bd} + q_{ac}\wedge q_{ac}.
\]
Now we turn to $\psi^1_{abcd}$.
\begin{multline*}
\psi^1_{abcd}=q_{cd}\wedge q_{cd} + h_{bc}^d\wedge h_{bc}^d + q_{cd}\wedge q_{cd}+h_{ac}^d\wedge h_{ac}^d + q_{ab}\wedge q_{ab} + h_{da}^b\wedge h_{da}^b+q_{ab}\wedge q_{ab}+h_{ca}^b\wedge h_{ca}^b\\
=q_{ad}\wedge q_{ad}+q_{bc}\wedge q_{bc} + k_{da}^{cb}\wedge k_{da}^{cb}+k_{cb}^{da}\wedge k_{cb}^{da}.
\end{multline*}
Therefore 
\[
\psi_{abcd}=4 k_{ab}^{cd}\wedge k_{ac}^{db}+k_{ad}^{bc}\wedge k_{ad}^{bc}+k_{bc}^{ad}\wedge k_{bc}^{ad}.
\]
We see that $\psi_{abcd}$ depends only on the angles and its sum over the manifold is zero. If $2$ is invertible in $A$ then
\[
\frac{1}{4}\psi_{abcd}=k_{ab}\wedge k_{ac}.
\]
This specializes to $z\wedge(1-z)$ in the case of an ideal hyperbolic tertrahedron from Remark \ref{rem:bloch}.

\bibliography{refs}

\end{document}